\documentclass[letterpaper, 10 pt, conference]{IEEEtran}  

\IEEEoverridecommandlockouts                              

\usepackage[normalem]{ulem}

\usepackage[latin1]{inputenc}
\usepackage[english]{babel}
\usepackage{amsfonts,amssymb,amsmath,comment}
\usepackage[dvipsnames]{xcolor}
\usepackage{enumerate}
\usepackage{tikz}
\usetikzlibrary{decorations, decorations.text,backgrounds}
\usepackage{graphicx}

\usepackage{cite}
\usepackage{url}
\usepackage{placeins}
\usepackage{color}
\usepackage{epsfig} 
\usepackage{times} 
\usepackage{cite}
\usepackage{algorithm}
\usepackage{algpseudocode}
\usepackage{bm}
\usepackage{epstopdf}

\usepackage{cite}

\newtheorem{theorem}{Theorem}[section]
\newtheorem{proposition}[theorem]{Proposition}

\newtheorem{definition}[theorem]{Definition} \newtheorem{lemma}[theorem]{Lemma}

\newtheorem{assumption}[theorem]{Assumption}
 
\title{\LARGE \bfseries
 Limited Rate Distributed Weight-Balancing and Average Consensus Over Digraphs\vspace{-0.3cm}
}

\author{Chang-Shen Lee, Nicol\`{o} Michelusi, and Gesualdo Scutari\vspace{-0.3cm}
\thanks{Lee and Michelusi are with the School of Electrical and Computer Engineering, Purdue University, West Lafayette, IN, USA. Scutari is with the  School of Industrial Engineering, Purdue University, West Lafayette, IN, USA. Emails:  {\tt\small
<lee2495,michelus,gscutari>@purdue.edu.} This work was supported by USA NSF under Grants CIF 1632599 and CIF 1719205; and in part by the ONR under Grant N00014-16-1-2244.}
}

\begin{document}

\maketitle
\thispagestyle{empty}
\pagestyle{empty}

\begin{abstract} 
Distributed quantized weight-balancing and average consensus over fixed digraphs are considered. A digraph with non-negative weights associated to its edges is weight-balanced if, for each node, the sum of the weights of its out-going edges is equal to that of its incoming edges. This paper proposes and analyzes the first distributed algorithm that solves the weight-balancing problem using only finite rate and simplex communications among nodes (compliant to the directed nature of the graph edges). Asymptotic convergence of the scheme is proved and a convergence rate analysis is provided. Building on this result, a novel distributed algorithm is proposed that solves the \emph{average} consensus problem over digraphs, using, at each iteration, finite rate simplex communications between adjacent nodes -- some bits for the weight-balancing problem, other for the average consensus. Convergence of the proposed quantized consensus algorithm to the average of the \emph{real} (i.e., unquantized) agent's initial values is proved, both almost surely and in $r$th mean for all positive integer $r$. Finally, numerical results validate our theoretical findings.
\end{abstract}
\section{Introduction}
Weight-balanced directed graphs -- digraphs wherein the sum of the  weights of the edges outgoing from each node is equal to the sum
of the weights of the edges incoming to the node -- play a key role in a number of network applications, including   distributed optimization \cite{Lorenzo2016}, distributed flow-balancing \cite{Hadjicostis2017}, distributed averaging and cooperative control \cite{Olfati-Saber2004}, just to name a few. 
In particular,  distributed average consensus over (di)graphs whereby  agents aim  at  agreeing on the sample average of their local values  has received considerable attention over the  years;   some applications include load-balancing \cite{Cybenko1989}, vehicle formation \cite{Fax2004}, and sensor networks \cite{Schizas2008,Scutari08}. 
Several of the aforementioned distributed algorithms, when run on digraphs,  
  require some form of graph regularity condition, such as the weight-balanced property (see, e.g., \cite{Rikos2014}). 

A variety of centralized algorithms have been proposed in the literature to balance a weighted digraph; see, e.g., \cite{Loh1970} and references therein. In this paper, we are interested in  the design of {\it distributed} iterative algorithms that solve the weight-balancing problem as well as the average consensus problem over digraphs, using only {\it quantized} information, {\it simplex communications}  among nodes (compliant to the directed nature of the graph edges), and without knowledge of the graph topology (with exception of the direct neighbors).
This is motivated by realistic scenarios where     (wired or wireless) communications on physical channels might not be full-duplex (e.g., nodes   transmit at different power and/or communication
channels are not symmetric due to interference) and are subject to rate constraints, meaning that only a finite number of bits can be reliably transmitted per channel use.   
To date, the design of such algorithms remains a challenging and open problem, as documented next. 

\subsection{Related works}
\noindent\textbf{Distributed weight-balancing}:
Distributed algorithms aimed at solving the weight-balancing problem were proposed in \cite{Gharesifard2012,Rikos2014,Hadjicostis2017,Rikos2016}. More specifically, \cite{Rikos2014} (resp. \cite{Gharesifard2012}) considered the real and integer (resp. discrete) case; in \cite{Hadjicostis2017}, the authors extended the real weight-balancing scheme of \cite{Rikos2014} to deal with box constraints on the graph weights. With the exception of \cite{Rikos2016}, all the aforementioned algorithms require communication with infinite rate. In fact, they transmit either real valued quantities or some unbounded integer information on the local balance\footnote{We use the term "balance" to denote the local imbalance with sign (positive or negative), and the term "imbalance" to denote its absolute value, cf. Definition \ref{def:eps}.}. While compliant with finite rate constraints, the distributed integer weight-balancing algorithm \cite{Rikos2016} requires {\it full-duplex} edge communications -- each agent must exchange information with \emph{both its out-neighbors and its in-neighbors} -- which may not comply with the underlying directed nature of the edges. Thus, to the best of our knowledge, algorithms that solve the weight-balancing problem using \emph{finite rate and simplex} communications are still missing.

\noindent\textbf{Distributed average consensus}:
Distributed average consensus algorithms have a long history, tracing back to the seminal works  \cite{Tsitsiklis1984,Xiao2004,Olfati-Saber2004}. All these early works assumed that agents can reliably exchange  unquantized information.   To cope with the limited data rate constraint, quantization was later introduced in consensus algorithms  and its effect  analyzed in  
\cite{Kashyap2007,Nedic2009,Aysal2008,Kar2010,Rajagopal2011,Li2011,Thanou2013,Wang2011,Zhu2015,Lee2017conf}, with \cite{Kashyap2007,Nedic2009,Aysal2008,Kar2010,Rajagopal2011,Li2011,Thanou2013} considering undirected graphs and \cite{Lee2017conf} directed but {\it balanced} digraphs only. Specific features of these algorithms are briefly discussed next. In \cite{Kashyap2007,Nedic2009}, 
 agents store and communicate quantized information; deterministic uniform quantization is adopted, so that   only converge to the average of the initial values of the agents' variables  {\it within some error} can be achieved. In \cite{Aysal2008}, agents utilize dithered (probabilistic) quantization to communicate with each other;  consensus   at one of the quantization values is achieved almost surely. In addition, the expected value of the consensus is equal to the average of the agents' initial data. Distributed quantized consensus algorithms converging to the   {\it exact} (i.e., unquantized) average of the  initial values of agents' variables  were proposed in \cite{Rajagopal2011,Li2011,Thanou2013}.  However, all these schemes are applicable only to {\it undirected} graphs.   Referring to the literature dealing with quantized  consensus over digraphs \cite{Lee2017conf,Wang2011,Zhu2015}, either quantization with infinite number of bits is considered \cite{Zhu2015} or weight-balancedness of the  digraph is needed \cite{Lee2017conf} to achieve the \emph{exact} average consensus (in contrast, \cite{Wang2011} does not converge to the exact average). Thus, to the best of our knowledge, algorithms that solve the \emph{exact} consensus problem using \emph{finite rate and simplex communications over unbalanced digraphs} are still missing.

\subsection{Summary of the main contribution}
The analysis of the literature shows that there are no  distributed algorithms  solving the weight-balancing and the {\it exact} average consensus problems (the latter over unbalanced digraphs), using {\it quantized information and simplex communications}. 
This paper provides an answer to these open questions. 

The first contribution is a novel distributed and quantized weight-balancing algorithm whereby agents transfer part of their balance -- 
the difference between the out-going and the incoming sum-weights, which should be zero for a weight-balanced graph -- to their neighbors via quantized signals, so as to reduce their own local imbalance -- the absolute value of their balance. 
We prove that the proposed scheme converges logarithmically to a weight-balanced solution. The developed convergence analysis is a novel technical contribution of the paper, and it is highlighted next. 
\begin{enumerate} 

\item First, we
prove that the total imbalance decreases iff. agents with positive balance transfer part of their balance to agents with negative one, termed the ``\emph{decreasing event}'' (cf. Lemma \ref{lemma_decr_event} and Definition \ref{def_event_D});
therefore, agents with positive balance closer to agents with negative balance are more \emph{important} than those farther away, since they more directly contribute to the decrement of the total imbalance.

\item Based on these findings, 
 we introduce a sophisticated metric (cf. (\ref{U_k})), function of the balances of agents, using the idea of decimal system representation, so that the balance of \emph{agents of higher importance} (i.e., closer to the nodes with negative balance) is represented by \emph{more significant digits}. Utilizing the proposed metric guarantees the occurrence of the ``\emph{decreasing event}'' within finite time, hence the decrement of the total imbalance (cf. Proposition \ref{lemma_minimal_decrease}).
\item  Second, we propose a novel diminishing step-size rule (cf. Assumption \ref{assump_gamma}), which guarantees that the balance at each node can be expressed as an integer multiple of the current step-size; 
  we show that this novel step-size design greatly facilitates the convergence analysis, since it allows one to tightly control the amount of decrement of the total imbalance at each stage.

\end{enumerate}

Building on the   above result, we then introduce a novel distributed algorithm that performs average consensus and weigh-balancing 
\emph{on the same time scale}  using only two-bit simplex communications -- one bit is devoted to the consensus and the other one to balance the digraph. Convergence of the agents' local variables to the {\it exact} average of the initial values is proved, both  in  mean square error sense and almost surely, along with  (deterministic) convergence of the sequence of weights to a   weight-balancing solution. 

The rest of the paper is organized as follows. In Section~\ref{sec:background}, we introduce some basic notation and preliminary definitions. Section \ref{section_weight_balance} introduces the proposed (one-bit)  quantized distributed weight-balancing algorithm along with its convergence properties. 
Section~\ref{section_consensus} presents the proposed  distributed two-bit quantized algorithm solving the average consensus problem while balancing the digraph, and study its convergence. Some numerical results are discussed in Section~\ref{section_simulation}, while  Section~\ref{section_conclusion} draws some conclusions. Due to space limitations, only sketches of the proofs are provided. The complete proof can be found in \cite{Lee_tr2018}.

\section{Notation and Background}\label{sec:background}
\subsection{Notation}
The set of real, integer, nonnegative integer, and postive integer numbers is denoted by $\mathbb{R}$, $\mathbb{Z}$, $\mathbb{Z}_+$, and $\mathbb{Z}_{++}$, respectively. The $0-1$ indicator function is denoted by $\mathcal{I}\{\mathcal{A}\}$:  the function returns 1 if the input argument $\mathcal{A}$ is true, and $0$ otherwise. We denote the probability space of a stochastic process $\{{\bf z}(k)\}_{k \in \mathbb{Z}_+}$ by $\left(\Omega, \sigma, P\right)$, where $\Omega$ is a sample space, $\sigma$ is a $\sigma$-algebra, and $P$ is a probability measure. In addition, filtration is denoted by $\mathcal{F} = \{\sigma(k)\}_{k \in \mathbb{Z}_+}$, where $\sigma(k)$ is a sub-$\sigma$-algebra of $\sigma$ for every $k \in \mathbb{Z}_+$. $\mathbb{E}[\cdot]$ denotes the expectation, the distribution with respect it is taken will be clear from the context. Vectors (resp. matrices) are denoted by lower-case (resp. capital), bold letters. Finally, all equalities and inequalities involving random variables are tacitly assumed to hold almost surely (i.e., with probability $1$), unless otherwise stated.
\subsection{Basic graph-related definitions} \label{def_graph}
Consider a network with $N$ agents, modeled  as a static, directed graph $\mathcal{G} = \{\mathcal{V},\mathcal{E}\}$, where $\mathcal{V} = \{1,\cdots,N\}$ is the set of vertices (the agents), and $\mathcal{E} \subset \mathcal{V} \times \mathcal{V}$ is the set of edges (the communication links). A directed edge from $i\in \mathcal{V}$ to $j\in\mathcal{V}$ is denoted by $(i,j) \in \mathcal{E}$, over which information flows. We assume that $\mathcal{G}$ does not contain self-loops, that is, $(i,i)\notin \mathcal{E}$. The  {\it in-neighbors} of node $i$ are nodes in  the set $\mathcal{N}_i^- = \{j:(j,i) \in \mathcal{E}\}$, while its out-neighbors are those in the set $\mathcal{N}_i^+ = \{j:(i,j) \in \mathcal{E}\}$. The cardinality of $\mathcal{N}_i^-$ (resp. $\mathcal{N}_i^+$) is called {\it in-degree} (resp. {\it out-degree}) of node $i$ and is denoted by $d_i^- = |\mathcal{N}_i^-|$ (resp. $d_i^+ = |\mathcal{N}_i^+|$). We denote  by $d(i,j)$ the directed distance between $i$ and $j\in\mathcal V$, that is, the length of the shortest path from $i$ to $j$; we set $d(i,i)=0$, for all $i\in \mathcal{V}$.  
We will consider  strongly connected digraphs. 
\begin{definition} \label{def:graph}  A digraph $\mathcal{G}$ is strongly connected if, for every two distinct nodes  $i,j\in \mathcal V$, there exists a directed path connecting $i$ to $j$, i.e. $d(i,j)<\infty,\forall i,j$.
\end{definition}  

Associated with the digraph $\mathcal{G}$, we define a weight matrix compliant to it, along with some related quantities instrumental to formulate the weight-balancing problem. 
\begin{definition}[Weight matrix]\label{def:MC} Given a digraph $\mathcal{G}$, a matrix   $\mathbf{A}\triangleq (a_{ij})_{i,j=1}^N\in \mathbb{R}^{N\times N}$ is said to be compliant to $\mathcal{G}$ if, for all $i,j\in\mathcal V$,  $$a_{ij}=\begin{cases}
\geq 0, & \text{if }(j,i)\in \mathcal{E};\\
0, & \text{otherwise}.
\end{cases}    $$
\vspace{-0.3cm}
\end{definition}
In the following, we will only consider compliant weight matrices.
\begin{definition}[In-flow and out-flow]  \label{def:S}
	Given a digraph $\mathcal{G}$ with weight matrix $\mathbf{A}$, the total {\it in}-flow and {\it out}-flow of node $i\in\mathcal V$ are defined as  $S_i^- \triangleq  \sum_{j \in \mathcal{N}_i^-}{a_{ij}}$ and $S_i^+ \triangleq \sum_{j \in \mathcal{N}_i^+}{a_{ji}}$, respectively. 
\end{definition}
 
 \begin{definition}[Node weight (im)balance] \label{def:eps}
 	Given a digraph $\mathcal{G}$ with weight matrix $\mathbf{A}$, the {\it weight balance} $b_i$ of node $i$ is defined as $b_i\triangleq S_i^--S_i^+$, and its {\it weight imbalance} as $\epsilon_i\triangleq |b_i|$.
	  The weight imbalance vector collecting the  $\epsilon_i$'s across the network  {is  $\bm{\epsilon}\triangleq (\epsilon_i)_{i=1}^N$}.
	\smallskip 
 \end{definition}
 
\begin{figure}[t]
	\centering
	\includegraphics[width = 1.2 in]{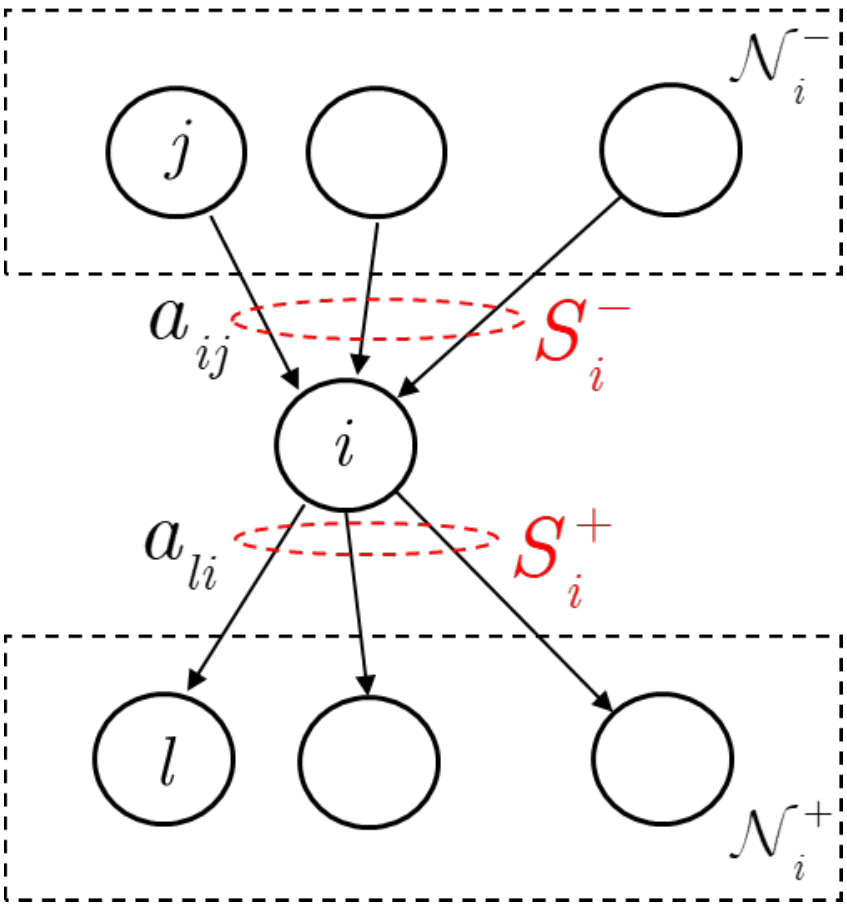}
	\caption{Some basic graph definitions.}		\label{system_model}%
	\vspace{-0.5cm}
\end{figure}

 \begin{definition}[Weight-balanced digraph]\label{def:WB}
 A digraph $\mathcal G$ is said to be \emph{weight-balanced} if its associated weight matrix $\mathbf A$ induces a total imbalance equal to zero, $\bm{\epsilon}=\mathbf{0}$.
 \end{definition}
Fig.~\ref{system_model} summarizes some of the quantities defined above. 

\section{Distributed One-Bit Weight-Balancing } \label{section_weight_balance}
In this section, we introduce a distributed, iterative algorithm to solve the weight-balancing problem using only quantized information and simplex communications. We are given a strongly connected digraph $\mathcal{G}$. Note that strong connectivity guarantees the existence of a matrix, compliant to the digraph  $\mathcal{G}$ (cf.~Definition~\ref{def:MC})  that makes $\mathcal{G}$ weight-balanced (cf.~Definition~\ref{def:WB}) \cite{Loh1970}.   Each node, say $i$, controls the set of weights $(a_{ij})_{j\in \mathcal{N}_i^-}$ associated with its incoming edges; the goal is to update iteratively the weights so that, eventually, they all converge to a matrix, compliant to $\mathcal{G}$,  which makes $\mathcal{G}$ weight-balanced. To do so, agents exchange information with their neighbors, under the following communication constraints: i) information flows according to the edge directions (simplex communications); and ii) information flows are quantized with a finite number of bits. We will denote by   
${\bf A}(k)=(a_{ij}(k))_{i,j=1}^N,  S_i^+(k) = \sum_{j \in \mathcal{N}_i^+}{a_{ji}(k)}, S_i^-(k) = \sum_{j \in \mathcal{N}_i^-}{a_{ij}(k)}, b_i(k) = S_i^-(k)-S_i^+(k), \text{ and } {\pmb\epsilon}(k) = (\epsilon_i(k))_{i=1}^N$ 
 the values of the associated variables at iteration $k$ of the algorithm (cf. Definitions \ref{def:MC}-\ref{def:eps}). We also denote the step-size at time $k$ as $\gamma(k)$. The proposed algorithm is formally stated in Algorithm~\ref{alg_dwb} and discussed next.
\begin{algorithm}[t]
\caption{Distributed Quantized Weight-Balancing} 
\label{alg_dwb} 
\begin{algorithmic} 
	 \State \textbf{Initialization}: Weight matrix ${\bf A}(0)=(a_{ij}(0))_{i,j=1}^N,$ with $a_{ij}(0)=1$ if $j\in \mathcal{N}_i^-$ and $a_{ij}(0)=0$ otherwise; step-size $\{\gamma(k)\}_{k \in \mathbb{Z}_+}$. 
	 \State Set $k=0$;
    \State \texttt{(S.1)} If ${\bf A}(k)$ satisfies a termination criterion: STOP;
    \State \texttt{(S.2)} Each agent $i$ broadcasts $n_i(k)$ to $\mathcal N_i^+$, where
    \begin{align}
    n_i(k) &= \mathcal{I}\{b_i(k) \geq d_i^+\cdot \gamma(k)\}. \label{signal_wb}
    \end{align}
    \State \texttt{(S.3)} Each agent $i$ collects signals $n_j(k)$ from $j\in \mathcal{N}_i^-$, and updates $\left(a_{ij}(k)\right)_{j \in \mathcal{N}_i^-}$ as    \begin{align}
    \!\!\!\!\!\!a_{ij}(k+1) &= a_{ij}(k)+n_j(k)\cdot\gamma(k), \quad  j \in \mathcal{N}_i^-. \label{a_ij_wb}\\
    \!\!\!\!\!\!\!b_i(k+1) &= b_i(k)-\gamma(k)\,d_i^+\,n_i(k)\!+\!\gamma(k)\!\!\sum_{j\in \mathcal{N}_i^-}\!\!n_j(k).\label{b_i_alg}
    \end{align}
\end{algorithmic}
\vspace{-0.3 cm}
\end{algorithm}
In \texttt{S.2}, each agent $i$ generates the binary signal $n_i(k)$
by comparing its weight balance $b_i(k)$ with the threshold $d_i^+\gamma(k)$; then, it broadcasts such signal to its out-neighbors. 
    In \texttt{S.3}, each agent $i$ collects the signals from its in-neighbors, and updates the corresponding weights according to \eqref{a_ij_wb}:   if  $n_j(k)>0$,   the   incoming weight $a_{ij}(k)$ will be increased by $\gamma(k)\,n_j(k)$. The balance of each agent is then updated according to (\ref{b_i_alg}). Roughly speaking, by \eqref{a_ij_wb}-(\ref{b_i_alg}) there is a {\it transfer} of the positive  balance among nodes in the network: 
    if $n_i(k)> 0$, the quantity $\gamma(k)\,d^+_i\,n_i(k)$ is subtracted from the balance $b_i(k)$ of node $i$ [cf.~(\ref{b_i_alg})], and equally divided among its out-neighbors $j\in\mathcal{N}_i^+$, which will increase their incoming weight   $a_{ji}(k)$ by $\gamma(k)\,n_i(k)$ [cf.~\eqref{a_ij_wb}]. 
    Note that Algorithm~\ref{alg_dwb}  is fully distributed: each agent $i$ only needs to know its out-degree $d_i^+$ and the initial value $b_i(0)=S_i^-(0)-S_i^+(0)$,  and to agree on a common step-size rule $\{\gamma(k)\}_{k \in \mathbb{Z}_+}$.

Before stating the convergence results, we first introduce the following step-size sequence $\{\gamma(k)\}_{k \in\mathbb{Z}_+}$.
\begin{assumption} \label{assump_gamma}
$\{\gamma(k)\}_{k \in\mathbb{Z}_+}$ is given by:
\begin{align*}
	\gamma(k) &= \left\{1,\frac{1}{2},\frac{1}{2},\frac{1}{4},\frac{1}{4},\frac{1}{4},\frac{1}{4},\frac{1}{8}, \cdots\right\} \nonumber \\
	 &= 2^{-n}, \text{ with } n\in \mathbb{Z}_+: 2^n-1 \leq k \leq 2^{n+1}-2.
\end{align*} 
\end{assumption}
The proposed step-size possesses three desired properties: (i) the balance at each node is always an integer multiple of the current step-size, which guarantees sufficient amount of decrease of the total imbalance $\Vert \pmb\epsilon(k)\Vert_1$, when decrease occurs (cf. Lemma \ref{lemma_decr_event}); (ii) it is diminishing, which prevents the algorithm from terminating prematurely; (iii) it is non-summable (cf. Lemma \ref{lemma_gamma_ub}), which together with the above two properties, guarantees $\Vert \pmb\epsilon(k)\Vert_1 \to 0$. Note that properties (ii) and (iii), i.e., $\lim_{k \rightarrow \infty}{\gamma(k)} = 0$ and $\sum_{k \geq 0}{\gamma(k)} = \infty$
are not surprising, since the one-bit signal (\ref{signal_wb}) can be regarded as a noisy version of the desired information $b_i(k)$. With such noisy information at hand, the use of a diminishing (nonsummable) step-size is consistent with similar choices adopted, e.g., in stochastic optimization \cite{Bertsekas1989}. However, typical choices in the context of stochastic optimization, e.g., $\gamma(k){=}1/(k{+}1), \forall k{\in}\mathbb{Z}_+$, do not satisfy (i); for these cases, it is hard to obtain theoretical performance guarantees for the weight-balancing problem, as opposed to the proposed step-size.

We are now ready to state our main convergence result. 

\begin{theorem}\label{thm_wb}
Let $\mathcal{G}$ be a strongly connected digraph. Let $\{\mathbf{A}(k)\}_{k\in \mathbb{Z}_+}$ be the sequence generated by  Algorithm \ref{alg_dwb}, with   step-size $\{\gamma(k)\}_{k \in \mathbb{Z}_+}$ satisfying Assumption \ref{assump_gamma}. Then, the following hold:

\noindent \texttt{(a) Asymptotic convergence}:
\begin{align}
\lim\limits_{k\rightarrow \infty}{{\bf A}(k)} = {\bf A}^\infty, \label{A_conv}
\end{align}
where ${\bf A}^\infty$  makes the digraph weight-balanced;

\noindent \texttt{(b) Convergence rate}:
\begin{align}
\Vert \pmb\epsilon(k)\Vert_1 = O\left(\frac{1}{k}\right). \label{eq_conv_rate}
\end{align}
\end{theorem}
\subsection{Proof of Theorem~\ref{thm_wb} (sketch)} \label{subsec_wb_proof}
\noindent {\bf Proof of statement (a):}  We prove the statement in two steps, namely: 1) 
 we show that  the total imbalance $\Vert \pmb\epsilon(k)\Vert_1$ is asymptotically vanishing; and 2) the sequence  $\{\mathbf{A}(k)\}_{k\in \mathbb{Z}_+}$ is convergent. Step~2 implies convergence   whereas Step~1 guarantees that the limit point of $\{\mathbf{A}(k)\}_{k\in \mathbb{Z}_+}$
is a solution of the weight-balancing problem. 

\noindent\underline{Step 1}:    $\lim\limits_{k \rightarrow \infty}{\Vert \pmb\epsilon(k)\Vert_1} = 0$. 
We begin  identifying the event $\mathcal D_k$ that ensure  $\Vert \pmb\epsilon(k)\Vert_1$ to strictly decrease.
\begin{definition}[Decreasing event $\mathcal D_k$] 
Let $\mathcal D_k$  be the ``\emph{decreasing event}'' defined as 
	\begin{equation}
	 \exists  i  \in \mathcal V \text{ and } j \in \mathcal{N}_i^+: 	 n_i(k)>0, \text{ and } b_j(k)<0.    	\label{event_D}\end{equation} \label{def_event_D}\vspace{-0.4cm}
\end{definition}
This event occurs when a node with sufficiently large \emph{positive} balance--node $i$ in \eqref{event_D}--triggers the update of the weights of an out-neighboring node with \emph{negative} balance--node $j$.
Indeed, we   show next that $\Vert \pmb\epsilon(k)\Vert_1$ decreases iff. $\mathcal D_k$ occurs, and remains unchanged otherwise.
\begin{lemma} \label{lemma_decr_event} There holds
 \begin{equation}\label{eq:lemma3}
 	\Vert{\pmb\epsilon}(k+1)\Vert_{1} \begin{cases}
\leq\Vert{\pmb\epsilon}(k)\Vert_{1}-2\,\gamma(k), & \text{if }\mathcal{I}\left(\mathcal D_k\right)=1,\\
=\Vert{\pmb\epsilon}(k)\Vert_{1}, & \text{otherwise.}
\end{cases}
 \end{equation}
\end{lemma}

Clearly, by Lemma \ref{lemma_decr_event}, we infer that $\Vert{\pmb\epsilon}(k)\Vert_{1}$ is not increasing.
However, this alone does not guarantee  $\Vert \pmb\epsilon(k)\Vert_1$ to vanish asymptotically; in fact, the decreasing event must occur sufficiently often. We  show next that, indeed,
the decreasing event occurs \emph{at least} once within a time window of finite duration (Proposition \ref{lemma_minimal_decrease2}).  
This together with the non-increasing and non-summable property of $\{\gamma(k)\}_{\mathbb{Z}_+}$,
will be enough to show that 
$\Vert \pmb\epsilon(k)\Vert_1$   asymptotically vanishes.

\begin{proposition} \label{lemma_minimal_decrease2}  
If $\Vert{\pmb \epsilon}(k)\Vert_1 \geq 2N(N-1)\gamma(k)$, then
$\mathcal I(\mathcal D_t)=1$ for some $k\leq t\leq k+N^{2N}$.
\end{proposition}
This proposition along with Lemma \ref{lemma_decr_event} and the decreasing nature of $\gamma(k)$ implies that 
\begin{proposition} \label{lemma_minimal_decrease}  
If $\Vert{\pmb \epsilon}(k)\Vert_1 \geq 2N(N-1)\gamma(k)$, then
 \begin{equation}\label{eq:error_decrease}
\left\Vert {\pmb \epsilon}\left(k+N^{2N}\right)\right\Vert_1 \leq \Vert {\pmb \epsilon}(k)\Vert_1-2\gamma\left(k+N^{2N}\right).
\end{equation}
\end{proposition}

\begin{IEEEproof}
We now prove Propositions \ref{lemma_minimal_decrease2} and \ref{lemma_minimal_decrease}.
The following lemma is instrumental to our proof.
\begin{lemma} \label{lemma_node_flow}
Let $\mathcal{V}^+(k) = \{i \in \mathcal{V}:b_i(k) \geq 0\}$ and $\mathcal{V}^{--}(k) = \{i \in \mathcal{V}:b_i(k) < 0\}$ denote the set of nodes with non-negative and negative balance at iteration $k$.  If $\mathcal{I}(\mathcal D_k)=0$,   then $\mathcal{V}^+(k+1) = \mathcal{V}^+(k)$, $\mathcal{V}^{--}(k+1) = \mathcal{V}^{--}(k)$.
\end{lemma}

\noindent In words, if the event $\mathcal D_k$ does not occur at iteration $k$, the sets of nodes having non-negative 
and those having negative balance do not change from $k\to k+1$. 
 
   Suppose that  $\Vert \pmb\epsilon(k_0)\Vert \geq 2N(N-1)\gamma(k_0)$, for some $k_0\in \mathbb{Z}_+$; let $T\triangleq \min \left\{k\in \mathbb{Z}_{+}\,:\,\right.$ $\left.\mathcal{I}\left(\mathcal{D}^{k_0+k}\right)=1\right\}$ be the (possibly, infinite) delay for the event $\mathcal D_k$ to occur for the first time since $k_0$. 
Invoking Lemma~\ref{lemma_decr_event}, we have \begin{align}
\Vert {\pmb \epsilon}(k_0+T+1)\Vert_1 \leq \Vert {\pmb \epsilon}(k_0)\Vert_1-2\gamma(k_0+T). \label{eps_dec}
\end{align}
 Suppose that  $T$ is bounded,   that is, $T\leq \bar{T}$, for some $\bar{T} < \infty$ (a fact that will be proved later). This means  that   $\mathcal D_k$ must occur at least once within  $\left[k_0,k_0+\bar{T}\right]$. We can write 
\begin{align}
&\!\!\!\!\!\!\Vert {\pmb \epsilon}\left(k_0+\bar{T}+1\right)\Vert_1 \stackrel{(\ref{eq:lemma3})}{\leq} 
\Vert {\pmb \epsilon}(k_0+T+1)\Vert_1 \nonumber \\
&\!\!\!\!\!\!\!\!\!\stackrel{(\ref{eps_dec})}{\leq} \!\!\Vert {\pmb \epsilon}(k_0)\Vert_1\!-\!2\gamma\!\!\left(k_0\!+\!T\right)
 \!\stackrel{(a)}{\leq}\! \Vert {\pmb \epsilon}(k_0)\Vert_1\!-\!2\gamma\left(k_0\!+\!\bar{T}\!+\!1\right),
\label{dec_lb}
\end{align}
where in $(a)$ we used the fact that  $\gamma(k)$ is non-increasing. 

It remains to prove that such a $\bar{T} < \infty$ exists, and in particular $\bar T=N^{2N}-1$, so that the proof of both propositions follows. Let $k\in \left[k_0,k_0+T-1\right]$. This implies that $\mathcal{V}^+(k)$ is invariant over the interval $\left[k_0,k_0+T\right]$ (Lemma~\ref{lemma_node_flow}); we thus write $\mathcal{V}^+ \triangleq  \mathcal{V}^+(k), \mathcal{V}^{++} \triangleq  \mathcal{V}^{++}(k)$ and $\mathcal{V}^{--} \triangleq  \mathcal{V}^{--}(k)$ for $k\in\left[k_0,k_0+T\right]$.
Clearly, we can partition the nodes in $\mathcal{V}^+$ based on their \emph{distance} to agents with negative balance. To this end, let
\begin{align}
\mathcal{V}_n \triangleq  \left\{i: i\in \mathcal{V}^+\text{ and } \underset{j \in \mathcal{V}^{--}}\min d(i,j) = n\right\},\ 
n\in\mathbb{Z}_+, \label{V_n}
\end{align} 
which represents the set (possibly empty) of agents  in $\mathcal{V}^+$ that are $n$-hops (directed) away from an agent with negative balance [recall that  $d(i,j)$ denotes the directed distance from $i$ to $j$, cf. Sec.~\ref{def_graph}].
Based on our intuitive discussion,
 given the total balance of agents in $\mathcal{V}^+$, the event $\mathcal D_k$ will occur sooner if the balance is concentrated on the agents in $\mathcal{V}_n$ with smaller $n$. 
In other words, nodes in $\mathcal V_1$ have a more direct impact on the occurrence of the decreasing event than nodes in $\mathcal V_2$; and so on,
  nodes in $\mathcal V_n$ have a more direct impact on the occurrence of the decreasing event than nodes in $\mathcal V_{n+1}$.
This observation implies that the distribution of the balance within $\mathcal{V}^+$ affects the time to the occurrence of $\mathcal D_k$. Therefore, 
our strategy to bound $T$ is to construct a metric $U(k)$, function of balance, 
representing how directly nodes in $\mathcal{V}^+$ influence the occurrence of the decreasing event.
Specifically, we design $U(k)$ for $[k_0, k_0+T]$, as a function of $\{b_i(k)\}_{i=1}^N$, with the following properties:
\begin{enumerate}[(a)]
\item $U(k)$ is strictly increasing;
\item The increments of $U(k)$ are integer numbers;
\item $U(k)$ is lower and upper bounded;
\end{enumerate}
Since $\gamma(k)$ can be regarded as the unit of balance, $b_i(k)/\gamma(k)$ can be regarded as the normalized balance. To satisfy (a), note that when agent $i$ in $\mathcal{V}_{n+1}$ triggers the update, some balance will be transferred from it to one of its out-neighbor, say, agent $j$ in $\mathcal{V}_n$. Hence, the increasing amount of $U(k)$ caused by the increased balance of agent $j$ should dominate the decreasing amount caused by the decreased balance of agent $i$. Since the normalized balance of agent $i$ will be decreased by $d_i$ and the normalized balance of agent $j$ will be increased by 1, this operation is similar to the ``carry'' operation in the decimal system: if a digit reaches $10$, then the next more significant digit will be increased by $1$. Based on the above observation, we propose a novel metric $U(k)$ to model the balance in $\mathcal{V}^+$ by leveraging the idea of "decimal system". Roughly speaking, $U(k)$ can be regarded as a "number" in which the total normalized balance of agents in $\mathcal{V}_1$ is represented by its most significant digit, the one of agents in $\mathcal{V}_2$ is represented by its second most significant digit, and so on. Since $U(k)$ is strictly increasing by integer quantities and upper bounded within $[k,k+T]$, it must follow that $T$ is finite. Formally, for $k\in\left[k_0,k_0+T\right]$
\begin{align}
U(k) &= \sum\limits_{n = 1}^{n_{\max}}{U_n\sum\limits_{i \in \mathcal{V}_n}{\min\left\{\frac{b_i(k)}{\gamma(k)},d_i^+\right\}}}, \text{ where} \label{U_k} \\
U_n &= \prod_{m = n+1}^{n_{\max}}{u_m}, u_m = 1+\sum\limits_{i \in \mathcal{V}_m}{d_i^+}, \label{U_n} \text{ and} \\
n_{\max} &= \max \left\{n:|\mathcal{V}_n| > 0\right\}. \nonumber
\end{align}
Unlike the decimal system, the base for the $m$th significant digit, $U_m$, is not fixed, as opposed to ten in the decimal system. The $m$th significant digit in $U(k)$, i.e., $\sum_{i \in \mathcal{V}_n}{\min\left\{b_i(k)/\gamma(k),d_i^+\right\}}$, is computed by first normalizing the balance of agents in $\mathcal{V}_n$ by $\gamma(k)$, which is a non-negative integer since $b_i(k)$ is a multiple of $\gamma(k)$, then clipping the excess part compared to $d_i^+$ and adding them together. Note that the clipping step is necessary to ensure that the resulting digit is no greater than its base.
The following lemma states three important properties of $U(k)$:
\begin{lemma} \label{lemma_Uk}
$U(k)$ exhibits the following properties:
\begin{enumerate}[(i)]
\item $U(k)$ is non-negative, $U(k)\geq 0$;
\item if $\Vert{\pmb \epsilon}(k)\Vert_1 \geq 2N(N-1)\gamma(k)$, $U(k)$ strictly is increasing, $U(k+1) \geq U(k)+1$;
\item $U(k)$ is upper-bounded by $U(k)< N^{2N}$.
\end{enumerate}
\end{lemma}
Invoking Lemma~\ref{lemma_Uk},  we readily get $N^{2N}\!>\!U(k_0+T)\geq U\left(k_0\right)+T \geq T$. Therefore, $T\leq \bar{T}=N^{2N}-1<\infty$.
This concludes the proof of Proposition~\ref{lemma_minimal_decrease}.

\begin{figure}[t!]
	\centering
	\includegraphics[width = \linewidth]{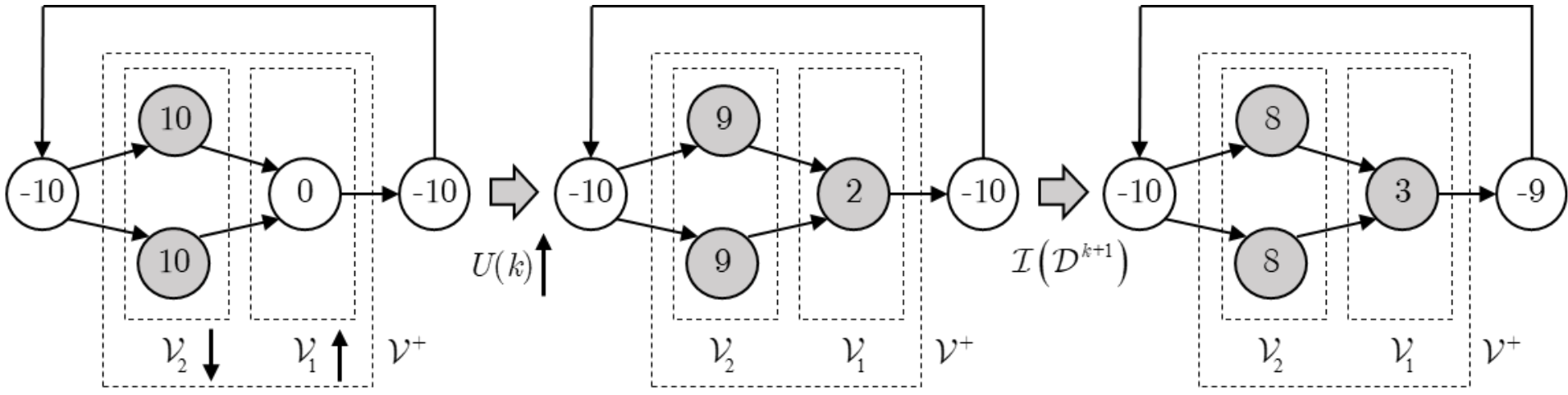}
	\vspace{-0.5cm}
	\caption{Illustration of the dynamics of $U(k)$ and $\Vert \pmb\epsilon(k)\Vert_1$. In this figure, the number in each node represents its balance, and $\gamma(k) = \gamma(k+1) = \gamma(k+2) = 1$; shaded nodes trigger the update. At time $k$, nodes belong to $\mathcal{V}_2$ trigger the updates, which causes some balance transferred from $\mathcal{V}_2$ to $\mathcal{V}_1$, making $U(k)$ increase. At time $k+1$, the increased metric $U(k+1)$ ensures the occurrance of $\mathcal{D}^{k+1}$, which causes $\Vert \pmb\epsilon(k+1)\Vert_1$ to decrease.}	\label{V_++}
	\vspace{-0.5cm}
\end{figure} 

Equipped with Lemma \ref{lemma_decr_event} and Proposition~\ref{lemma_minimal_decrease}, we can now prove that $\Vert\pmb\epsilon(k)\Vert_1\to 0$.

Since $\Vert\pmb\epsilon(k)\Vert_1$ is a non-increasing sequence, it is sufficient to show that
$\forall \epsilon>0$ there exists $k_\epsilon\in \mathbb{Z}_+$ such that
 $\Vert\pmb\epsilon(k_\epsilon)\Vert_1<\epsilon$.
To this end, let $k_0\in \mathbb{Z}_+$ be such that $2N(N-1)\gamma(k_0)\leq\epsilon$; note that such $k_0$ exists, as
$\lim_{k\to\infty}\gamma(k)=0$. If $\Vert\pmb\epsilon(k_0)\Vert_1<2N(N-1)\gamma(k_0)$,  the result follows readily, with $k_\epsilon=k_0$. Now suppose that $\Vert\pmb\epsilon(k_0)\Vert_1\geq 2N(N-1)\gamma(k_0)$. 
Since $\gamma(k)\leq\gamma(k_0)$ for $k\geq k_0$, it suffices to show that
$\exists k_{\epsilon}>k_0:\Vert\pmb\epsilon(k_\epsilon)\Vert_1<2N(N-1)\gamma(k_\epsilon)$.
We  prove it by contradiction.  Suppose  that $\Vert\pmb\epsilon(k)\Vert_1\geq2N(N-1)\gamma(k)$, for all $k\geq k_0$.
  Applying recursively \eqref{eq:error_decrease} yields 
$\Vert\pmb\epsilon\left(k_0+t\,N^{2N}\right)\Vert_1
 \leq \Vert\pmb\epsilon(k_0)\Vert_1-2\sum_{\tau=1}^t\gamma\left(k_0+\tau N^{2N}\right), \forall t\in \mathbb{Z}_+$. Taking the limit $t\to\infty$, yields $0 \leq \Vert\pmb\epsilon(k_0)\Vert_1-2 \sum_{\tau=1}^\infty\!\gamma\left(k_0+\tau N^{2N}\right)=-\infty$.
\end{IEEEproof}

\noindent\underline{Step 2: } We show that $\{{\bf A}(k)\}_{k \in \mathbb{Z}_+}$ is convergent.
Since it is nondecreasing [cf.~(\ref{a_ij_wb})], 
 it suffices to prove that $\{{\bf A}(k)\}_{k \in \mathbb{Z}_+}$  is bounded, which is proved in the following lemma.
\begin{lemma} \label{lemma_a_bound}
The sequence  $\{{\bf A}(k)\}_{k \in \mathbb{Z}_+}$ is bounded. 
\end{lemma}

\noindent {\bf Proof of statement (b):} 
We want to show that
\begin{align} 
\exists 0\!<\!M\!<\!\infty \text{ and }\,\bar{k}\in \mathbb{Z}_{++}\!: k\Vert\pmb\epsilon(k)\Vert_1\!<\!M, \forall k\geq \bar{k}. \label{limit2}
\end{align}
    The following two lemmas provide some intermediate results that will be invoked to prove \eqref{limit2}.
\begin{lemma} \label{lemma_gamma_ub}
$\frac{1}{k+1}\leq\gamma(k) \leq \frac{2}{k+1}, \forall \in \mathbb{Z}_+$.
\end{lemma}

\begin{lemma} \label{lemma_T-}
For every ${k} \in \mathbb{Z}_+$, there exists $k_1 \geq  {k}$ such that $\Vert \pmb\epsilon(k_1)\Vert < 2N(N-1)\gamma(k_1)$.
\end{lemma}

We show next that  \eqref{limit2} holds with the following choice of $M$ and $\bar{k}\in \mathbb{Z}_{++}$:     
\begin{align*}
M= 4N(N-1)\max\left\{1,T_0\right\}, \text{ with }\tilde{T}\triangleq3N(N-1)N^{2N},
\end{align*}
and $\bar{k}$ satisfying
\begin{subequations}
\begin{align}
&\!\!\!\!\!\!\!\bar{k}\geq k_0 + \tilde{T}, \label{k_bar1}\\
&\!\!\!\!\!\!\!\Vert \pmb\epsilon(t)\Vert_1 \!<\! 2N(N\!-\!1)\gamma(t), \text{ for some integer } t\in[k_0,\bar{k}],\label{k_bar2}
\end{align}
\end{subequations}
where  $k_0 \in \mathbb{Z}_{++}$ is such that  \begin{align}
\gamma(k)=\gamma(k_0),\ \forall k\in[k_0,k_0+\tilde{T}]. \label{large_k_2}
\end{align}
Note that such a $k_0\in \mathbb{Z}_{++}$ and $\bar{k}\in \mathbb{Z}_{++}$  exist, since the time interval during which $\gamma(k) = \frac{1}{2^n}$ has duration
 $\left[\left(2^{n+1}-1\right)-1\right]-\left[\left(2^{n}-1\right)\right]+1 =2^n$, tends to infinity as $n \to \infty$ (and thus $k\to \infty$). This implies that one can always find a sufficiently large $k_0$ such that (\ref{large_k_2}) is satisfied. A similar argument along with Lemma~\ref{lemma_T-} can be used to show that \eqref{k_bar2} also holds, for sufficiently large  ${k}_0$ and  $\bar{k}$. Hence, such a $\bar{k}$ exists.
Let
\begin{align}
\!\!\!\!\tilde{T}_k\triangleq\min\{t\in[0,k]\!:\!\Vert \pmb\epsilon(k\!-\!t) \Vert_1 \!<\! 2N(N\!-\!1)\gamma(k\!-\!t)\}. \label{T}
\end{align}
The existence of such $\tilde T_k$ is guaranteed by \eqref{k_bar2}.

We show next that $\tilde{T}_k$ is upper bounded by $\tilde{T}$, i.e.,  $\tilde{T}_k\leq \tilde{T}$. This holds trivially if $\tilde{T}_k=0$, hence let us consider the case $\tilde{T}_k>0$. By definition of $\tilde{T}_k$ and $\bar{k}$ [cf.~\eqref{k_bar2}] there hold:
\begin{align}
	&\!\!\!\!\!\!\Vert \pmb\epsilon(k\!-\!\tilde{T}_k) \Vert_1 \!<\! 2N(N\!-\!1)\gamma(k\!-\!\tilde{T}_k),\text{ with }  k-\tilde{T}_k\geq k_0,\label{eq:contrad_0}\\
	&\!\!\!\!\!\!\Vert \pmb\epsilon(k\!-\!\tilde{T}_k\!+\!t) \Vert_1 \!\geq\! 2N(N\!-\!1)\gamma(k\!-\!\tilde{T}_k\!+\!t), \forall t\in [1,\,\tilde{T}_k].\label{eq:contrad_1} 
\end{align}
Note that $k-\tilde{T}_k\geq k_0$ follows from  \eqref{k_bar2}.
We will now prove $\tilde{T}_k\leq \tilde{T}$ via contradiction.
    Suppose that $\tilde{T}_k>\tilde{T}$, and consider the following lemma.
    \begin{lemma} \label{lemma_bound_intv_large_k} 
Let $\tau\geq k_0$ be such that 
    $\Vert \pmb\epsilon(\tau-1)\Vert_1 < 2N(N-1)\gamma(\tau-1)$ and $\Vert \pmb\epsilon(\tau)\Vert_1 \geq 2N(N-1)\gamma(\tau)$.
    Then, there exists an integer 
    $t\in  [1, \tilde{T}]$   such that $\Vert \pmb\epsilon(\tau+t)\Vert_1 < 2N(N-1)\gamma(\tau+t)$.
\end{lemma}

It is clear that (\ref{eq:contrad_0})-(\ref{eq:contrad_1})
satisfy the  conditions of the lemma with $\tau=k-\tilde{T}_k$.
Therefore,  $\Vert \pmb\epsilon(k-\tilde{T}_k+t)\Vert_1 < 2N(N-1)\gamma(k-\tilde{T}_k+t)$ for some $t\in  [1, \tilde{T}]$,
which  contradicts (\ref{eq:contrad_1}).

We have thus proved that $\tilde{T}_k\leq \tilde{T}$.
Using this upper bound, we can write 
\begin{align}\label{eq:upper_error_case_2}
\Vert {\pmb \epsilon}(k)\Vert_1
\stackrel{(\ref{eq:lemma3})}{\leq}
\Vert {\pmb \epsilon}(k-\tilde{T}_k)\Vert_1
 \stackrel{(a)}{\leq}  
 2N(N-1)\frac{2}{k-\tilde{T}_k+1},
\end{align}
where $(a)$ comes from \eqref{T} and Lemma \ref{lemma_gamma_ub}. Note that $k-\tilde{T}_k+1>0$. Using  $\tilde{T}_k\leq \tilde{T}$  in \eqref{eq:upper_error_case_2}, yields  
\begin{align}
k\cdot\Vert {\pmb \epsilon}(k)\Vert_1
&\leq
 2N(N-1)\frac{2k}{k-\tilde{T}+1}\leq M, \forall k\geq \bar{k},
\end{align}
which completes the proof. \hfill $\square$

\section{Distributed Quantized Average Consensus} \label{section_consensus}
In this section,  we devise  a  novel  distributed algorithm  solving 
the quantized average consensus problem over (non-balanced) digraphs; the proposed scheme  builds  on the distributed quantized weight-balancing algorithm (Algorithm~\ref{alg_dwb}) introduced in Section~\ref{section_weight_balance}, as described next.      

Consider the same network setting as in Section~\ref{section_weight_balance}. Let
$y_i(0)\in\mathbb R$ denote the initial sample owned by agent $i$. The goal is to design a distributed algorithm whereby agents will eventually agree on the average   of the initial values,
\begin{align}\label{eq:consensus}
\bar{y}(0) \triangleq \frac{1}{N}\sum_{i = 1}^N{y_i(0)}.
\end{align}
   Agents can exchange quantized information with their neighbors via simplex communications. Since the digraph is not assumed to be balanced, plain consensus schemes (using quantization) cannot be readily used; a weight-balancing procedure needs to be incorporated in the consensus updates. The proposed idea is then to combine the weight-balancing algorithm introduced in Section~\ref{section_weight_balance} with the average  consensus protocol based on probabilistic quantization, which we recently proposed in \cite{Lee2017conf}. 
   The   new algorithm is designed so that these two building blocks  run \emph{on the same time-scale}. 
    The   scheme   is summarized in Algorithm~\ref{alg_dac} and works as follows. 
    
    Every agent $i$ owns two sets of variables, namely: the weights associated to the in-neighbors $(a_{ij})_{j \in \mathcal{N}_i^-}$ and the local  estimate $y_i$ aiming at asymptotically converging to (\ref{eq:consensus}); we denote by $(a_{ij}(k))_{j \in \mathcal{N}_i^-}$ and  $y_i(k)$ the value of these variables at time $k$. At each iteration  $k$, based upon its current balance $b_i(k)$ and local estimate $y_i(k)$, agent $i$ generates and broadcasts to its out-neighbors the quantized signals $n_i(k)$ and $x_i(k)$ (Step~2). More specifically, the signal $n_i(k)$ is generated according to \eqref{signal_wb_alg_ac} while, to generate $x_i(k)$, agent $i$ first clips its local estimate $y_i(k)$ within the quantization range $[q_{\min},q_{\max}]$ [cf.~(\ref{clippedY})], and then quantizes the clipped estimate $\tilde{y}_i(k)$ via \eqref{quantizer} to build  $x_i(k)$.
 Upon receiving the signals $(n_j(k))_{j \in \mathcal{N}_i^-}$ and $(x_j(k))_{j \in \mathcal{N}_i^-}$ from its in-neighbors, agent $i$ updates its weights $(a_{ij}(k))_{j \in \mathcal{N}_i^-}$ using the quantized weight-balancing rule  introduced in Algorithm~\ref{alg_dwb} [cf. \eqref{a_ij}], and the local variable $y_i(k)$ according to \eqref{update_consensus}. The update in \eqref{update_consensus} aims at forcing a consensus among  the local variables $y_i(k)$ on the average $\bar{y}(0)$. In fact, the third term in \eqref{update_consensus} is instrumental to align the local copies $y_i(k)$ while the second term $ b_i(k)x_i(k)$ is a correction needed to preserve the average of the iterates, i.e., $(1/N)\sum_i y_i(k+1)=(1/N)\sum_i y_i(k)$, for all $k\in \mathbb{Z}_+$, which guarantees that,   if all the  $y_i(k)$ are asymptotically  consensual, it must be $\big|y_i(k)-(1/N)\sum_i y_i(k)\big|=\big|y_i(k)-(1/N)\sum_i y_i(0)\big| \underset{k\to\infty}{\longrightarrow} 0$.  
   Note that all the above  steps in the algorithm can be implemented in a distributed fashion, using only local information.
\begin{algorithm}[h]
\caption{Distributed Quantized Average Consensus and Weight-Balancing} 
\label{alg_dac} 
\begin{algorithmic} 
	\State \textbf{Data}: ${\bf A}(0)$ with $a_{ij}(0)=1$, if $j\in \mathcal{N}_i^-$, and $a_{ij}(0)=0$, otherwise; $\{\gamma(k)\}_{k \in \mathbb{Z}_+}, \{\alpha(k)\}_{k \in \mathbb{Z}_+}, q_{\min}<q_{\max}, {\bf y}(0)$. \\
	Set $k=0$;
    \State \texttt{(S.1)} If a termination criterion is satisfied: STOP;
    \State \texttt{(S.2)} Each agent $i$ broadcasts the following signals to its out-neighbors 
    \begin{align}
    n_i(k) &= \mathcal{I}\{b_i(k) \geq d_i^+\gamma(k)\}, \label{signal_wb_alg_ac} \\
    x_{i}(k) &=\left\{\begin{array}{*{20}l}
	q_{\max}, & \text{w.p. } p_i(k), \\
	q_{\min}, & \text{w.p. } 1-p_i(k), 
	\end{array}\right. \label{quantizer} \\
\text{where }p_i(k) &= \frac{\tilde{y}_i(k)-q_{\min}}{q_{\max}-q_{\min}}, \\
\tilde{y}_i(k) &= \min\bigr\{\max\bigr\{y_{i}(k),q_{\min}\bigr\},q_{\max}\bigr\}.\label{clippedY}
	\end{align}
    \State \texttt{(S.3)} Agent $i$ collects signals from its in-neighbors, and updates $\{a_{ij}(k)\}_{j \in \mathcal{N}_i^-},y_i(k)$:
        \begin{align}
    a_{ij}(k+1) &= a_{i,j}(k)+n_j(k)\gamma(k), \label{a_ij} \\
    y_i(k+1) &= y_i(k)
     +\alpha(k)\sum\limits_{j \in \mathcal{N}_i^-}{a_{ij}(k)\big(x_j(k)-x_i(k)\big)}    \nonumber \\
    &
    +\alpha(k)b_i(k)x_i(k).  \label{update_consensus}
    \end{align}
\end{algorithmic}
\vspace{-0.2 cm}
\end{algorithm}

We next introduce the  assumption on $\bar{y}(0)$ and the step-size used in the consensus updates.
\begin{assumption}[Informative $\bar{y}(0)$] \label{assump_sbar}
The average $\bar{y}(0)$ [cf.~(\ref{eq:consensus})]  
satisfies $\bar{y}(0) \in [q_{\min}, q_{\max}]$.
\end{assumption}
\begin{assumption} \label{assump_alpha_ms} The step-size 
   $\{\alpha(k)\}_{k \in \mathbb{Z}_+}$ satisfies:     
\begin{align*}
&\alpha(k)>0, \alpha(k+1) \leq \alpha(k), \forall k \in\mathbb{Z}_+, \\
&\sum\limits_{k=1}^{\infty}{\alpha(k)} = \infty, \sum\limits_{k=1}^{\infty}{\alpha(k)^2}< \infty.
\end{align*}
\end{assumption}

 It is important to remark that Assumption~\ref{assump_sbar} does not require each local data $y_i(0)$ to be confined within  the quantization range, nor does it require the range of $y_i(0)$ to be known. This is a major departure from the literature, which requires $y_i(0)$ to be within the quantization range -- see, e.g. \cite{Rajagopal2011,Li2011,Wang2011,Thanou2013}. We  instead need that the   \emph{average} $\bar y(0)$   falls within the quantization interval $[q_{\min}, q_{\max}]$, which is a less restrictive condition. 
 For example, if agents are estimating a common unknown parameter $\theta$, throughout  the measurement $y_i(0) = \theta+\omega_i$, where $\omega_i$ is zero mean Gaussian noise, i.i.d. across agents, then $\bar{y}(0)$ is the sample mean estimate across the agents. In this case, a bound on $y_i(0)$ is hard to obtain (theoretically it is unbounded), but the bound of the parameter, $\theta\in[\theta_{\min}, \theta_{\max}]$, is known in many cases.
Even worse, $\max_{i\in\{1,2,\dots,N\}}|y_i(0)|\to\infty$ for $N\to\infty$, whereas the sample average $\bar y(0)\to \theta$, so that the sample average $\bar y(0)$ becomes more and more informative for large $N$, whereas the initial local measurements become more and more unbounded. In this case, agents can thus simply set $(q_{\min},q_{\max}) = (\theta_{\min},\theta_{\max})$.
Herein, we are not interested in non-informative ${\bar y}(0)$. In fact, in this case $\bar{y}(0)$ does not provide information for estimating $\theta$.

We now state the convergence result of  Algorithm~\ref{alg_dac}.

\begin{theorem} \label{thm_ac}
Let $\mathcal{G}$ be a strongly connected digraph. Let $\left\{\mathbf{y}(k)=(y_i(k))_{i=1}^N\right\}_{k\in \mathbb{Z}_+}$ be the sequence generated by  Algorithm~\ref{alg_dac} under Assumptions \ref{assump_gamma}, \ref{assump_sbar}, and \ref{assump_alpha_ms}. Then, \\
\noindent \texttt{(a) Almost sure convergence}: 
\begin{align}
\mathbb{P} \left(\lim\limits_{k \rightarrow \infty}{{\bf y}(k) = \bar{y}(0)\cdot {\bf 1}} \right) =1. \label{eq_as_main}
\end{align}
\noindent \texttt{(b) Convergence in the mean square sense}:
\begin{align}
\lim\limits_{k \rightarrow \infty}{\mathbb{E}\left[\Vert{\bf y}(k) - \bar{y}(0) \cdot {\bf 1}\Vert^2\right]} = 0. \label{eq_ms_main}
\end{align}
\end{theorem}

\subsection{Proof of Theorem \ref{thm_ac} (sketch)} \label{subsec_ac_proof}
We begin by rewriting the dynamics of $y_i(k)$ and $x_i(k)$ in vector form.
Let ${\bf y}(k) = \left(y_i(k)\right)_{i=1}^N$, ${\bf x}(k) = \left(x_i(k)\right)_{i=1}^N$
 and $\tilde{\bf y}(k) \triangleq  \left(\tilde{y}_i(k)\right)_{i=1}^N$.
  Using \eqref{update_consensus},
we can express ${\bf y}(k+1)$ as
\begin{align}
{\bf y}(k+1)={\bf y}(k)-\alpha(k){\bf L}^{+}(k)\mathbf x(k),
\label{update_y}
\end{align}
where ${\bf L}^+(k) \triangleq  {\bf S}^+(k)-{\bf A}(k)$, with ${\bf S}^+(k) \triangleq {\rm diag}\left\{S_1^+(k), \cdots, S_N^+(k)\right\}$ and $S_i^+(k) = \sum_{j \in \mathcal{N}_i^+}{a_{ji}(k)}$.

Note that, due to the dithered quantization, $\mathbb E[\mathbf x(k)|\mathbf y(k)]=\tilde{\mathbf y}(k)$.
In order to investigate the dynamics of the consensus error, we introduce the following quantity.
\begin{align}
&V({\bf y}(k))\triangleq\Vert{\bf y}(k)-\bar{y}(0){\bf 1}\Vert^2. \label{V}
\end{align} 
  $V({\bf y}(k))$ satisfies the following dynamics.  
 \begin{lemma} \label{lemma_ac_error_relation}
In the setting  of Theorem~\ref{thm_ac}, there holds
\begin{align}
\!\mathbb E\left[V({\bf y}(k\!+\!1))|{\bf y}(k)\right] 
&\!= V({\bf y}(k)) -2\alpha(k){\bf y}(k)^T{\bf L}^+(k)\tilde{\bf y}(k) \nonumber\\
&\!+\!\alpha(k)^2{\mathbb E}\!\left[ \Vert {\bf L}^+(k){\bf x}(k)\Vert^2|{\bf y}(k)\right]\!. \label{eq: descent_I}
\end{align}
\end{lemma}
We bound below the second term on the RHS of the inequality (\ref{eq: descent_I}). This is instrumental to show that the negative term in (\ref{eq: descent_I}) is dominant with respect to the last term. 

\begin{lemma} \label{lemma_ac_dec_term}
There holds: 
\begin{align}\label{eq:lower_bound_negative_term}
{\bf y}(k)^T{\bf L}^+(k)\tilde{\bf y}(k) \geq c_1 V\left({\bf y}(k)\right) - c_2\Vert \pmb\epsilon(k)\Vert_1,
\end{align}
for some finite constants $c_1, c_2 > 0$. 
\end{lemma}
Using  (\ref{eq:lower_bound_negative_term}) in (\ref{eq: descent_I}) yields,
\begin{align}
\mathbb E\left[V({\bf y}(k\!+\!1))|{\bf y}(k)\right] \leq 
V({\bf y}(k))\!-\!2c_1\alpha(k)\left[V({\bf y}(k))\!-\!c(k)\right], \label{err_descent}
\end{align}
where in the last inequality we defined $c(k)\triangleq
\frac{c_2}{c_1}\Vert \pmb\epsilon(k)\Vert_1
+\alpha(k)\frac{c_3}{2c_1}$, for some deterministic constant $c_3 \geq \Vert {\bf L}^+(k)\Vert^2{\mathbb E}\left[ \Vert {\bf x}(k)\Vert^2|\mathbf y(k) \right]> 0$, and the last inequality follows from
the fact that ${\bf A}(k)$, hence ${\bf L}^+(k)$, are bounded
(cf. Lemma \ref{lemma_a_bound}), and $\mathbf x(k)$ is the output of a finite rate quantizer, hence it is bounded as well.

\noindent \textbf{Proof of statement (a):}
It is sufficient to show that  $V$ defined in {\eqref{V}} satisfies the conditions of \cite[Theorem 1]{Kar2010}, namely:   
\begin{align*}
&\underset{\left\Vert{\bf y}(k)-\bar{y}(0){\bf 1}\right\Vert \geq \epsilon}\inf V({\bf y}(k)) > 0, \forall \epsilon > 0, \\
&V(\bar{y}(0)\cdot{\bf 1}) = 0, \text{ and }
\underset{{\bf y}(k) \rightarrow \bar{y}(0)\cdot{\bf 1}}\limsup V({\bf y}(k)) = 0.\end{align*}
In addition, from \eqref{err_descent} we have that 
\begin{align*}
\mathbb E\left[V({\bf y}(k+1))|{\bf y}(k)\right]-V({\bf y}(k))\leq  g(k)-2c_1\alpha(k)V({\bf y}(k)), 
\end{align*}
where $g(k) = 2c_1\alpha(k)c(k)$ satisfies $\sum_{k \geq 0}{g(k)}< \infty$ and  $g(k) > 0$, since $\sum_{k\geq 0}{\alpha(k)^2}<\infty$ and $c(k)^2=O(1/k^2)$ (cf. Theorem \ref{thm_wb}). So far, we have shown that all conditions of \cite[Theorem 1]{Kar2010} are satisfied. Hence, statement (a) holds. 
\noindent \textbf{Proof of statement (b):} 
To prove statement (b), we use the following lemma:
\begin{lemma}\label{lemma_y_global_bound}
Let $\{{\bf y}(k)\}_{k \in \mathbb{Z}_+}$ be the sequence generated by Algorithm \ref{alg_dac}, in the setting of Theorem \ref{thm_ac}. Then, $y_{i,\min} \leq y_i(k) \leq y_{i,\max}$ for some $y_{i,\min}, y_{i,\max} \in \mathbb{R}$. 
\end{lemma}
Since $\left\vert y_i(k)-\bar{y}(0)\right\vert$ is bounded for all $i \in \mathcal{V}$, we have $\mathbb E\left[\Vert{\bf y}(k)-\bar{y}(0){\bf 1}\Vert^2\right] < \infty$, together with the almost sure convergence proved in statement (a) implies statement (b).

\section{Numerical Results} \label{section_simulation}
In this section, we present some numerical results to validate our theoretical findings. Agent's initial data $\left(y_i\right)_{i \in \mathcal{V}}$ are generated i.i.d.  according to the uniform distribution on $(q_{\min},q_{\max})$. The digraph modeling the network is generated as follows. There are  $N = 6$ nodes. A directed ring is first constructed linking all the nodes, so that the digraph is ensured to be  strongly connected. Then a directed edge on each pair of nodes is randomly added, with probability $p$; several values of $p$ are considered.  Roughly speaking, $p$ can be regarded as a sparsity measure: the graph is sparse if $p$ is small and dense if $p$ is large.  

We measure the performance of the algorithm  using the total imbalance, i.e., $\Vert \pmb\epsilon(k) \Vert_1$, to monitor the weight-balancedness whereas the consensus disagreement is measured in terms of  MSE, defined as
\begin{align*}
{\rm MSE} = \frac{1}{N}\sum\limits_{i = 1}^N{(y_i(k)-\bar{y})^2}.
\end{align*}
The simulation results are averaged over 100 realizations.

\begin{figure}[t] 
	\centering
	\vspace{-0.5cm}
	\includegraphics[width = \linewidth]{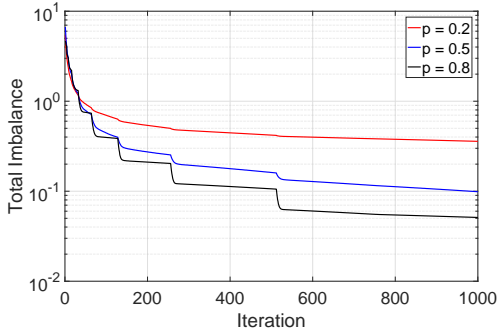} 
	\caption{Quantized weight-balancing problem: $\Vert \pmb\epsilon(k)\Vert_1$ by Algorithm~\ref{alg_dwb} vs. number of iterations. Different curve represents different graph sparsity.}	 \label{fig_imb}	
	\vspace{-0.2cm}
\end{figure}

In Fig.~\ref{fig_imb} we plot the total imbalance  $\Vert \pmb\epsilon(k)\Vert_1$ generated by Algorithm~\ref{alg_dwb} versus the number of iterations. Different curves refer to different level of sparsity of the graph (probability values $p$). The curves are averaged over $100$ independent graph realizations, with $100$ independent initial value realizations being evaluated in each graph realization. The following comments are in order. The total imbalance $\Vert \pmb\epsilon(k)\Vert_1$ is a non-increasing function of the iterations, which is consistent with our analytical results. The results  show that, as expected, the proposed algorithm performs better over denser graphs, since on denser graphs there are more communications and more frequent updates among agents. In addition, one can see that the curve of $\Vert \pmb\epsilon(k) \Vert_1$ can be partitioned into nearly flat line segments and steep line segments, for all cases. In the former case, $\Vert \pmb\epsilon(k) \Vert_1$ does not have large decrease since balance is mainly transferred among agents with positive balance; on the other hand, in the latter case, $\Vert \pmb\epsilon(k) \Vert_1$ has large decreases, since balance is mainly transferred from agents with positive balance and those with negative balance. 

\begin{figure}[t] 
	\centering
	\vspace{-0.2cm}
	\includegraphics[width = \linewidth]{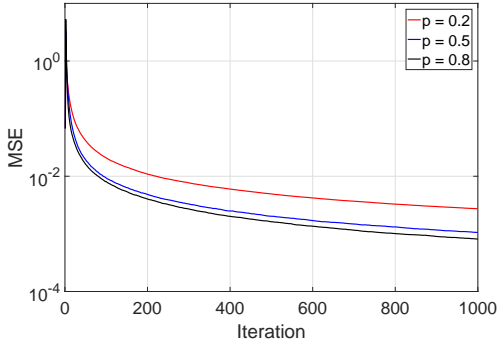}
	\caption{Quantized consensus problem: MSE generated   by Algorithm~\ref{alg_dac} vs. number of iterations. Different curve represents different graph sparsity.}		\label{fig_mse}
	\vspace{-0.5cm}
\end{figure}

Fig. \ref{fig_mse} shows the corresponding MSE performance, generated by Algorithm~\ref{alg_dac}. As expected, the denser graphs have better MSE performance,  due to i)  more frequent communications of the variables $\left(y_i\right)_{i \in \mathcal{V}}$ among agents; and ii)  more balanced graphs they experience (cf.~Fig.~\ref{fig_imb}).

\section{Conclusions} \label{section_conclusion}
In this paper, we  introduced a novel distributed algorithm that solves the weight-balancing problem. The proposed scheme uses quantized information (one-bit) and simplex communications. Asymptotic convergence was proved along with the convergence rate. Building on this result, a second contribution of the paper, was a novel distributed average consensus algorithm  over (non-balanced) digraphs that uses only two-bit simplex communications. Convergence of the algorithm was proved using a novel line of analysis: a metric inspired by the decimal system as well as a dedicated step-size are proposed to show that the total imbalance will eventually converge to zero. Finally, numerical results showed that the proposed algorithms perform well in practice.

\bibliographystyle{IEEEtran}
\bibliography{IEEEabrv,Weight_balancing_arxiv_091718}	
\addtolength{\textheight}{-12cm}   



\end{document}